\documentclass[twoside, 12pt]{article}
\usepackage{amssymb, amsmath, mathrsfs, amsthm}%, wasysym}
\usepackage{graphicx}
\usepackage{color}
\usepackage{float}

\newcommand{\bd}{\begin{description}}
\newcommand{\ed}{\end{description}}
\newcommand{\bi}{\begin{itemize}}
\newcommand{\ei}{\end{itemize}}
\newcommand{\be}{\begin{enumerate}}
\newcommand{\ee}{\end{enumerate}}
\newcommand{\beq}{\begin{equation}}
\newcommand{\eeq}{\end{equation}}
\newcommand{\beqs}{\begin{eqnarray*}}
\newcommand{\eeqs}{\end{eqnarray*}}

\definecolor{DarkGreen}{rgb}{0.2, 0.6, 0.3}

\newtheorem{theorem}{Theorem}
\newtheorem{conjecture}{Conjecture}

\newtheorem{lemma}{Lemma}

\newtheorem{claim}{Claim}

\newtheorem{problem}{Problem}

\begin{document}
\title{\textbf{Note on the connectivity keeping spiders in
$k$-connected graphs} \footnote{Supported by the National
Science Foundation of China (No. 12061059),
Doctoral research project of Tianjin Normal University (52XB2111)
and the Qinghai Key Laboratory of Internet of Things Project
(2017-ZJ-Y21).} }
\author{Meng Ji\footnote{Corresponding author: College of Mathematical Science, Tianjin Normal University, Tianjin, China {\tt
jimengecho@163.com}}, \ \ Yaping Mao\footnote{Faculty of Environment and Information
Sciences, Yokohama National University, 79-2 Tokiwadai, Hodogaya-ku,
Yokohama 240-8501, Japan. {\tt mao-yaping-ht@ynu.ac.jp}}}

\date{}
\maketitle

\begin{abstract}
W. Mader [J. Graph Theory 65 (2010), 61--69] conjectured that for
any tree $T$ of order $m$, every $k$-connected graph $G$ with
$\delta(G)\geq\lfloor\frac{3k}{2}\rfloor+m-1$ contains a tree
$T'\cong T$ such that $G-V(T')$ remains $k$-connected. In 2010,
Mader confirmed the conjecture for the $k$-connected graph if $T$ is
a path; very recently, Liu et al. confirmed the conjecture if
$k=2,3$. The conjecture is open for $k\geq 4$
till now. In this paper, we show that Mader's conjecture is true for
the $k+1$-connected graph if $T$ is a spider and $\Delta(G)=|G|-1$.
\\[2mm]
\end{abstract}

\section{Introduction}

In this paper, all the graphs are finite, undirected simple. For
graph-theoretical terminology and notations not defined here, we
follow \cite{BM}. Some basic symbols and definitions ones are needed
to be introduced. The minimum degree and the connectivity number of
a graph $G$ is denoted by $\delta(G)$ and $\kappa(G)$, respectively.
the vertex set and edge set of $G$ is denoted by $V(G)$ and $E(G)$,
respectively. For a subgraph $H\subseteq G$ and a subset
$V'\subseteq V(G)$, write $V'\cap H$ for $V'\cap V(H)$. For any
vertex $x\in G$, denote the set of neighbors of $x$ in $G$ by
$N_{G}(x)$ and for a subgraph $H\subseteq G$, define
$N_{G}(H):=\bigcup_{x\in H}N_{G}(x)-V(H)$. For $H\subseteq G$,
$G(H)$ stands for the subgraph induced by $H$ in $G$ and $G-H:=G-V(H)$. For
$H\subseteq G$, set $\delta_{G}(H):=\min_{x\in H}{d_{G}(x)}$.
For any edge $xy\in E(G)$ and $H\subseteq G$, write $H\cup xy$ for the
subgraph $(V(H)\cup \{x,y\},E(H)\cup{xy})$ of $G$. Let $x,y$-path be
a path $P$ from $x$ to $y$. For $u,v\in P$, let $P[u,v]=P[v,u]$ be
the subpath of $P$ between $u$ and $v$, and $P[u,v)$ means
$=P[u,v]-v$. Denotes the complete graph on vertex set $S$ by $K(S)$.

Let $G$ be a graph. A separating set $S\subseteq G$ denotes a vertex
cut-set, and the cardinality of a minimum separating set $S$ is
denoted by $|S|=\kappa(G)$. If $S$ is a minimum separating set of
the graph $G$, then the union components $F$ of $G-S$ with
$G-S-F\neq \emptyset$ is called a fragment $F$ to $S$, and the
complementary fragment $G-(S\cup V(F))$ is denoted by $\bar{F}$. If
a fragment of $G$ does not contain any other fragments of $G$, then
it is called an end of $G$. Clearly, every graph contains an end
except for the complete graphs. For a fragment $F$ of $G$ to $S$, it
definitely follows $S=N_{G}(F)$. For $S\subseteq G$, denote the
complete graph induced by $S$ by $K(S)$ and the graph $G[S]:=G\cup K(S)$.

In 1972, Chartrand, Kaigars and Lick \cite{CKL} proved that
\emph{every $k$-connected graph $G$ with $\delta (G)\geq
\lfloor\frac{3k}{2}\rfloor$ has a vertex $x$ with $\kappa(G-x)\geq
k$}. Nearly forty years later, Fujita and Kawarabayashi \cite{SK}
extended the result, that is, every $k$-connected graph $G$ with
$\delta (G)\geq \lfloor\frac{3k}{2}\rfloor+2$ has an edge $xy$ such
that $G-\{x,y\}$ remains $k$-connected. Furthermore, Fujita and
Kawarabayashi posed the following conjecture.

\begin{conjecture}[Fujita and Kawarabayash \cite{SK}]\label{fujita}
For all positive integers $k$, $m$, there is a (least) non-negative
integer $f_{k}(m)$ such that every $k$-connected graph $G$ with
$\delta(G)\geq \lfloor\frac{3k}{2}\rfloor+ f_{k}(m)-1$ contains a
connected subgraph $W$ of order $m$ such that $G-V(W)$ is still
$k$-connected.
\end{conjecture}

In 2010, Mader \cite{mader1} proved Conjecture \ref{fujita} is true
if $f_{k}(m)=m$ and $W$ is a path. Mader \cite{mader1} proposed the
following conjecture. The author will show it by using the Mader's
method.

\begin{conjecture}[Mader \cite{mader1}]\label{conj}
For every positive integer $k$ and finite tree $T$ of order $m$,
every $k$-connected finite graph $G$ with minimum degree
$\delta(G)\geq\lfloor\frac{3k}{2}\rfloor+|T|-1$ contains a subgraph
$T'\cong T$ such that $G-V(T')$ remains $k$-connected.
\end{conjecture}

In the case that $G$ is a graph with small connectivity, in 2009 Diwan and Tholiya \cite{AN} confirmed the conjecture for $k=1$; Tian et. al\cite{Tian1} and Tian et.al \cite{Tian2} studied some special trees, such as stars, double stars and path-star and so on, for $2$-connected graphs; Hasunuma \cite{H} proved some results on $2$-connected graphs with girth conditions; Hasunuma, Ono \cite{TK} and Lu, Zhang \cite{Lu} showed that it holds for $2$-connected graphs if $T$ is a tree with diameter condition; Very
recently, Hong and Liu \cite{HL2} showed that it holds for $2$-connected graphs if $T$ is a caterpillar or a spider. Furthermore,  Hong et. al\cite{HL} confirmed the conjecture for $k=2,3$.
In the case that $G$ is a graph with general connectivity, Mader \cite{mader1}
proved the conjecture is true if $T$ is a path; Mader \cite{mader2} showed that
the conjecture is true for $k$-connected graphs if $\delta\geq 2(k-1+m)^2+m-1$.
In this paper, the authors will confirm Mader's conjecture for the spider tree using Mader's structural theorem in Section $2$; in Section $3$ we conclude the paper.
In order to show our theorems, we need some structural lemmas
to prove our results.
%For fragments $F_{i}$ of $G$ to $S_{i}$ for $i=1, 2$, we define $S(F_{1},F_{2}) :=S_{1}\cap F_{2}\cup S_{1}\cap S_{2}\cup S_{2}\cap F_{1}$, where $\cup$ ties more than $\cap$.

\begin{lemma}[\cite{mader1}]\label{le1}
Let $G$ be a graph with $\kappa(G)=k$ and let $F$ be a fragment of
$G$ to $S$. Then it follows that if $F$ is an end of $G$ with $F\geq
2$, then $G[S]-V(\bar{F})$ is $(k+1)$-connected.
%(b) If $F_{1}\cap F_{2}\neq \emptyset$ holds, then $S(F_{1},F_{2})\geq k$ follows.\\
%(c) If $S(F_{1},F_{2})\geq k$, then $|F_{1}\cap S_{2}|\geq |\bar{F_{2}}\cap S_{1}|$ holds.
\end{lemma}

\begin{lemma}[\cite{mader1}]\label{le2}
Let $G$ be a $k$-connected graph and let $S$ be a separating set of
$G$ with $|S|=k$. Then the following holds.

$(a)$ For every fragment $F$ of $G$ to $S$, $G[S]-V(F)$ is
$k$-connected.

$(b)$ Assume $\delta(G)\geq \lfloor\frac{3k}{2}\rfloor +m-1$ and let
$F$ be a fragment of $G$ to $S$. If $W\subseteq G-(S\cup V(F))$ has
order at most $m$ and $\kappa(G[S]-V(F\cup W))\geq k$ holds, then
also $\kappa(G-V(W))\geq k$ holds.

$(c)$ Assume $(G,C)\in \mathcal{F}_{k}(m)$ and let $F$ be a fragment
of $G$ to $S$ with $C\subseteq G(F\cup S)$. If $W\subseteq G-(S\cup
V(F))$ has order at most $m$ and $\kappa(G[S]-V(F\cup W))\geq k$
holds, then also $\kappa(G-V(W))\geq k$ holds.
\end{lemma}

\begin{lemma}[\cite{mader1}]\label{le3}
For all $(G,C)\in \mathcal{F}_{k}^{+}(m)$ and $p_{0}\in G-V(C)$,
there is a path $P\subseteq G-V(C)$ of length $m-1$ starting from
$p_{0}$, such that $\kappa(G-V(P))\geq k$ holds.
\end{lemma}

\begin{theorem}[\cite{AN}]\label{1-connected}
Let $G$ be a connected graph with minimum degree $d\geq 1$. Then for
any tree $T$ of order $d$, $G$ contains a subtree $T'$ isomorphic to
$T$ such that $G-T'$ is connected.
\end{theorem}

%Mader \cite{mader1} showed that the minimum degree condition is tight for the Mader conjecture. However, we could relax the minimum condition and show it holds for a path. Let $G$ be a connected graph, then the degree sum condition of $G$ is defined as follows:
%$$\sigma_{2}(G) = \min\{d(x) + d(y)|x, y \in V(G), xy \notin E(G)\}.$$

\begin{table}[h]
\caption{Main known results on Mader conjuecture\,.} \label{tab:1} \centering
\begin{tabular}{cccc}
\cline{1-4} \hline $\kappa(G)$ & $\delta$  & $T$ & Authors and References\\[0.1cm]
\hline  $1$ & $\delta\geq m$ & $T$ is any tree & Diwan and Tholiya \cite{AN}; \\[0.1cm]
$2\ or\ 3$ & $\delta\geq m+2\ or\ m+3$ & $T$ is any tree & Hong and Liu \cite{HL}; \\[0.1cm]
$k$ & $\delta\geq \lfloor\frac{3k}{2}\rfloor+m-1$ & $T$ is a path &  Mader \cite{mader1}; \\[0.1cm]
$k$ & $2(k-1+m)^{2}+m-1$ & $T$ is any tree & Mader \cite{mader2};\\[0.1cm]
%$k$ & $\delta\geq \lfloor\frac{3k}{2}\rfloor+m-1$ & $T$ is a spider;\\[0.1cm]
\cline{1-4}
\end{tabular}
\end{table}

\section{Main results}
We define the set $\mathcal{F}_{k}(m)$ containing all pairs $(G,C)$ satisfying the following conditions:
\begin{itemize}
  \item $G$ is a $k$-connected graph with $|G|\geq k+1$;
  \item $C\subseteq G$ is a complete subgraph with $|C|=k$ and $\delta_{G}(G-V(C))\geq \lfloor\frac{3k}{2}\rfloor+m-1$;
  \item we denote by $\mathcal{F}_{k}^{+}(m)$ all the pairs $(G,C)\in \mathcal{F}_{k}(m)$ with $\kappa(G)\geq k+1$.

\end{itemize}

The spiders are considered and defined now. For a tree, if there is at most one vertex of degree at least $3$, then this tree is called a \emph{spider} (Specially, a path is also a spider; but we a star is not a spider). Each \emph{leg} of a spider is a path from the vertex adjacent to the root $x_{0}$ to a vertex of degree $1$; if there are $z$ legs, then denote the spider by $T^{t_{1},t_{2},\cdots,t_{z}}_{m}$, where $|T^{t_{1},t_{2},\cdots,t_{z}}_{m}|=m$ and $t_{i}$ denotes the order of the $i$th leg with $t_{1}+t_{2}+\cdots+ t_{z}+1=m$. If there are $t$ legs of order one, then we abbreviate $T^{1,1,\cdots,1,m-t-1}_{m}$ as $T^{t;m-t-1}_{m}$.

\begin{lemma}\label{lem1}
Let $t\geq 0$ be an integer. For any $(G,C)\in
\mathcal{F}_{k}^{+}(m)$ and any $s_{0}\in G-V(C)$ with $d(s_{0})=|G|-1$, there is a spider
$T^{t;m-t-1}_{m}\subseteq G-V(C)$ of order $m$ rooted at
$s_{0}$ such that $$\kappa(G-V(T^{t;m-t-1}_{m}))\geq k.$$
\end{lemma}

\begin{proof}
We perform an induction on
the order $n$ of the graph $G$ for the lemma. Clearly, the order of
the graph $G$ must no less than $\lfloor\frac{3k}{2}\rfloor+m$ since
$\delta(G-C)\geq \lfloor\frac{3k}{2}\rfloor+m-1$. Then it holds for
$(G,C)\in \mathcal{F}_{k}^{+}(m)$ if $G$ is a complete graph with
order at least $\lfloor\frac{3k}{2}\rfloor+m$. So we just need to
consider the case that $G$ is not complete and $|G|\geq
\lfloor\frac{3k}{2}\rfloor+m+1$. Now assume that $G$ is a graph with
smallest order and with $\Delta(G)=|G|-1$ such that $(G,C)\in \mathcal{F}_{k}^{+}(m)$ is a
counterexample to Lemma \ref{lem1} for $k$, and some $C\subseteq G$,
and $m$.

Subject to above assumption, we find out on
the order $m$ of tree $T^{t;m-t-1}_{m}$ satisfying above assumption. From Lemma \ref{le3}, there exists a path $P\subseteq G-V(C)$ of
length $m-1$ starting from $s_{0}$ such that $\kappa(G-V(P))\geq k$
holds. Let $P=\{s_{0}p_{1}\}\cup P[p_{1},p]$. Since $|N_{G}(u)\cap (G-P)|\geq \lfloor\frac{3k}{2}\rfloor+m-1-(|P|-1)\geq k$ for any vertex $u\in V(P[p_{1},p])$, then $\kappa(G-s_{0}p_{1})\geq k$, where $s_{0}p_{1}\in E(G)$
is a subpath of $P$ and also a spider $T_{2}^{1,0}$. Hence, Lemma \ref{lem1} holds when $m=1,2$. Now suppose that a maximal spider $T^{t;j}_{t+j+1}$ with root
$s_{0}$ and legs $s_{1},s_{2},\cdots, s_{t}$ satisfies the
following conditions.
\begin{itemize}
\item[(i)]$2\leq|T^{t;j}_{t+j+1}|=t+j+1<m$;

\item[(ii)]$\kappa(G-T^{t;j}_{t+j+1})\geq k$.
\end{itemize}

Note that $s_{1},s_{2},\cdots, s_{t}\in V(G)$ and $s_0s_i\in E(G)$,
$1\leq i\leq t$, and $T^{t;j}_{t+j+1}-\{s_{i}\,|\,1\leq
i\leq t\}$ is a path of order $j+1$, say
$P=s_{0}p_{1}\cup P[p_{1},p_{j}]:=p_{1}p_{2}\cdots p_{j}$. Simply, we set
$H=G-T^{t;j}_{t+j+1}$.
\begin{claim}\label{claim1}
$\kappa(H)=k$.
\end{claim}
\begin{proof}
Assume, to the contrary, that $\kappa(H)>k$. Since
$$
|N_{G}(x)\cap (H-C)|\geq \lfloor
3k/2\rfloor+m-1-k-(m-2)=\lfloor k/2\rfloor+1
$$
for $x\in \{s_{0},p_{j}\}$, it follows that there exists a vertex $s\in H-C$ such that
$xs\in E(G)$. Note that $T^{t;j}_{t+j+1}\cup
xs$ is a spider rooted at $s_{0}$ of order
$t+j+2\leq m$, which contradicts the choice
of $T^{t;j}_{t+j+1}$. Complete the proof of Claim \ref{claim1}.
\end{proof}

Since $H$ is not a complete graph, it follows that $|H|\geq k+2$. An
end $E$ is contained in $H$ with $E\cap C=\emptyset$. Set
$S:=N_{H}(E)$. Then $|S|=k$. Furthermore, let $\bar{E}=H-S-E$.

\begin{claim}\label{claim2}
$|E|\geq 2$.
\end{claim}
\begin{proof}
Assume, to the contrary, that $|E|=1$. It satisfies $k=d_{H}(z)\geq
\lfloor\frac{3k}{2}\rfloor+m-1-|T^{t;j}_{t+j+1}|\geq
\lfloor\frac{3k}{2}\rfloor$ for each $z\in E$, which means $k=1$ and
$|T^{t;j}_{t+j+1}|=m-1$. We get
$V(T^{t;m-t-2}_{m-1})\subseteq N_{G}(z)$ because of
$\delta(G)\geq \lfloor\frac{3k}{2}\rfloor+m-1$. Then
$|T^{t;m-t-1}_{m}|:=|T^{t;m-t-2}_{m-1}\cup
zx|=m$ for $x\in\{s_{0},p_{m-t-2}\}$ and $G-T^{t;m-t-1}_{m}=H-z$ is $1$-connected. And also,
$|T^{t+1;m-t-2}_{m}|:=|T^{t;m-t-2}_{m-1}\cup
zx|=m$ and $G-T^{t+1;m-t-2}_{m}=H-z$ is $1$-connected, which contradicts that
$T^{t;j}_{a}$ is a maximal spider. Complete the proof of Claim \ref{claim2}.
\end{proof}

From Claim \ref{claim2}, we have $|E|\geq 2$. Then the graph
$H[S]-\bar{E}$ is $(k+1)$-connected from Lemma \ref{le1}. From above
assumption, we know $\kappa(G)>k=\kappa(H)$, thus
$N_{G}(T^{t;j}_{t+j+1})\cap E\neq \emptyset$. Let $y$ be
one of the farthest vertices to $s_{0}$ on $T^{t;j}_{t+j+1}$
with $N_{G}(y)\cap E \neq \emptyset$. Otherwise, $S$ is also a separating set of $G$, which contradicts $\kappa(G)>k$. Suppose that $q$ is one vertex in $N_{G}(y)\cap E$. We
distinguish the following two cases to show this lemma. We will construct two lager spiders $T^{t;j+1}_{t+j+2}$ and $T^{t+1;j}_{t+j+2}$ such that $G-T^{t;j+1}_{t+j+2}$ and $G-T^{t+1;j}_{t+j+2}$ remains $k$-connected, respectively.

\begin{claim}\label{cl3}
There is a lager spider $T^{t;j+1}_{t+j+2}$ such that $G-T^{t;j+1}_{t+j+2}$ remains $k$-connected.
\end{claim}
\begin{proof}
Suppose that $y\in \{p_{1},p_{2},\cdots,p_{j},s_{0}\}$.
Let $\bar{P}=P[p_{j},y)$. Consider the graph $G-
(T^{t;j}_{t+j+1}-\bar{P}):=H \cup \bar{P}$. Since
$|N_{G}(x)\cap H|\geq \lfloor\frac{3k}{2}\rfloor+m-1-(t+j)\geq
\lfloor\frac{3k}{2}\rfloor+1\geq k$ for any $x\in V(\bar{P})$, it
follows that
$$
\kappa(G-(T^{t;j}_{t+j+1}-\bar{P}))=\kappa(H\cup
\bar{P}))\geq k.
$$
As $y$ is the farthest vertex to $s_{0}$ on
$T^{t;j}_{t+j+1}$, we have $N_{G}(\bar{P})\cap
E=\emptyset$. Naturally, $S$ is also a minimum separating set of $H
\cup \bar{P}$, and $E$ is an end of $H \cup \bar{P}$. From Lemma
\ref{le1}, $(H \cup \bar{P})[S](E\cup S)=H[S]-\bar{E}$ is
$(k+1)$-connected. Furthermore, it follows that $C\subseteq H\cup
\bar{P}-E$ and
$$
\delta_{(H\cup \bar{P})}(H\cup \bar{P}-C)\geq
\left\lfloor\frac{3k}{2}\right\rfloor+m-1-(t+j+1-|\bar{P}|)\geq
\left\lfloor\frac{3k}{2}\right\rfloor+|\bar{P}|.
$$
Hence, $(H\cup \bar{P},C)\in \mathcal{F}_{k}(|\bar{P}|+1)$ and
$(H\cup \bar{P}[S](E\cup S),K(S))\in
\mathcal{F}^{+}_{k}(|\bar{P}|+1)$. From Lemma \ref{le3}, there
exists a path $Q\subseteq E$ of order $|\bar{P}|+1$ starting from
$q$ such that $(H\cup \bar{P})[S](E\cup S)-Q$ is $k$-connected. Now
for $(H\cup \bar{P},C)\in \mathcal{F}_{k}(|\bar{P}|+1)$ by Lemma
\ref{le2}(c), then $\kappa((H\cup \bar{P})[S](E\cup S)-Q)\geq k$.
The spider
$T^{t;j+1}_{t+j+2}:=(T^{t;j}_{t+j+1}\setminus
\bar{P})\cup yq\cup Q$ rooted at $s_{0}$ has order
$|T^{t;j}_{t+j+1}|+1\leq m$ and
$G-V(T^{t;j+1}_{t+j+2})=H-Q$ is $k$-connected, a contradiction.

Suppose that $y\in \{s_{1},s_{2},\cdots,s_{t}\}$.
Let $\bar{P}=P[p_{j},p_{1})$ and then there is a spider
$T^{t;j+1}_{t+j+2}:=(T^{t;j}_{t+j+1}-\bar{P})\cup
yp\cup Q$ with order $t+j+2\leq m$ according to
the same way as above case.
\end{proof}

\begin{claim}\label{cl4}
There is a lager spider $T^{t+1;j}_{t+j+2}$ such that $G-T^{t+1;j}_{t+j+2}$ remains $k$-connected.
\end{claim}
\begin{proof}

Since $d_{G}(s_{0})=|G|-1$, then $N(s_{0})\cap E\neq \emptyset$. Let $N(s_{0})\cap E=\{q\}$. Furthermore, from Lemma \ref{le3}, there
exists a vertex $q\subseteq E$ such that $H[S](E\cup S)-q$ is $k$-connected. Hence, $G-V(T^{t+1;j}_{t+j+2})=H-q$ is $k$-connected, a contradiction.
\end{proof}
By Claim \ref{cl3} and Claim \ref{cl4}, we completed the proof.
\end{proof}

\begin{theorem}\label{thm1}
Every $k+1$-connected graph $G$ with
$\delta(G)\geq\lfloor\frac{3k}{2}\rfloor+m-1$ and $\Delta(G)=|G|-1$ for positive integers
$k,m,t$, contains a spider $T^{t;m-t-1}_{m}$ such that
$$\kappa(G-V(T^{t;m-t-1}_{m}))\geq k.$$
\end{theorem}
\begin{proof}
Clearly, the complete graph with order at least $\lfloor\frac{3k}{2}\rfloor+m$ holds. Suppose that $G$ is $k+1$-connected with $\delta(G)\geq \lfloor\frac{3k}{2}\rfloor+m-1$
and not complete. 
%Then there is an end $E$ in $G$ and let $S=N_{G}(E)$ with $|S|=k$.

%\setcounter{fact}{0}
%\begin{fact}\label{fact3}
%If $\kappa(G)=k$, then $(G[S]-\bar{E},K(S))\in \mathcal{F}^{+}_{k}(m)$.
%\end{fact}
%\begin{proof}
%If $|E|=1$, then it follows from Theorem \ref{1-connected} that
%$k=d_{G}(x)\geq \lfloor\frac{3k}{2}\rfloor+m-1\geq
%\lfloor\frac{3k}{2}\rfloor$ for each $x\in E$ implies $k=m=1$. By
%Theorem \ref{1-connected}, there exists one vertex, say $v\in V(G)$,
%such that $G-v$ is $1$-connected, and hence $(G[S]-\bar{E},K(S))\in
%\mathcal{F}^{+}_{1}(1)$. Suppose $|E|\geq 2$. Since $G[S]-\bar{E}$
%is $(k+1)$-connected by Lemma \ref{le1} and $\delta_{G}(E)\geq
%\lfloor\frac{3k}{2}\rfloor+m-1$, it follows that
%$(G[S]-\bar{E},K(S))\in \mathcal{F}^{+}_{k}(m)$.
%\end{proof}

%By Fact \ref{fact3} and Lemma \ref{lem1} there exists a spider
%$T^{t;m-t-1}_{m}\subseteq E$ of order $m$ such that $G[S]-\bar{E}-
%T^{t;m-t-1}_{m}$ remains $k$-connected. From Lemma \ref{le2}$(b)$,
%we have $\kappa(G-T^{t;m-t-1}_{m})\geq k$.

\setcounter{case}{1}
Assume that $\kappa(G)\geq k+1$.
We again find a maximal spider $T^{t;j}_{t+j+1}$ with root
$s_{0}$ and legs $s_{1},s_{2},\cdots, s_{t}$ satisfying the
conditions: $(i)$ $2\leq|T^{t;j}_{t+j+1}|=t+j+1<m$; $(ii)$ $\kappa(G-T^{t;j}_{t+j+1})\geq k$.
 Then in the following proof we prove exactly in the same way and symbols as
in the proof of Lemma \ref{lem1}. Naturally, we can show
$$\kappa(H)=k \ and \ |E|\geq 2.$$

Since $H:=G-T^{t;j}_{t+j+1}$ is not a complete graph, it follows that $|H|\geq k+2$. An
end $E$ is contained in $H$ with $E\cap C=\emptyset$. Set
$S:=N_{H}(E)$. Then $|S|=k$. Furthermore, let $\bar{E}=H-S-E$.

Then the graph
$H[S]-\bar{E}$ is $(k+1)$-connected from Lemma \ref{le1}. From above
assumption, we know $\kappa(G)>k=\kappa(H)$, thus
$N_{G}(T^{t;j}_{t+j+1})\cap E\neq \emptyset$. Let $y$ be
one of farthest vertices to $s_{0}$ on $T^{t;j}_{t+j+1}$
with $N_{G}(y)\cap E \neq \emptyset$. Set $q\in N_{G}(y)\cap E$.

\setcounter{claim}{0}
\begin{claim}\label{cl5}
There is a lager spider $T^{t;j+1}_{t+j+2}$ such that $G-T^{t;j+1}_{t+j+2}$ remains $k$-connected.
\end{claim}
Suppose that $y\in \{p_{1},p_{2},\cdots,p_{j},s_{0}\}$.
Let $\bar{P}=P[p_{j},y)$. Consider the graph $G-
(T^{t;j}_{t+j+1}-\bar{P}):=H \cup \bar{P}$. Since
$|N_{G}(x)\cap H|\geq \lfloor\frac{3k}{2}\rfloor+m-1-(t+j)\geq
\lfloor\frac{3k}{2}\rfloor+1\geq k$ for any $x\in V(\bar{P})$, it
follows that
$$
\kappa(G-(T^{t;j}_{t+j+1}-\bar{P}))=\kappa(H\cup
\bar{P}))\geq k.
$$
As $y$ is the farthest vertex to $s_{0}$ on
$T^{t;j}_{t+j+1}$, we have $N_{G}(\bar{P})\cap
E=\emptyset$. Naturally, $S$ is also a minimum separating set of $H
\cup \bar{P}$, and $E$ is an end of $H \cup \bar{P}$. From Lemma
\ref{le1}, $(H \cup \bar{P})[S](E\cup S)=H[S]-\bar{E}$ is
$(k+1)$-connected. Furthermore, it follows that $C\subseteq H\cup
\bar{P}-E$ and
$$
\delta_{(H\cup \bar{P})}(H\cup \bar{P}-C)\geq
\left\lfloor\frac{3k}{2}\right\rfloor+m-1-(t+j+1-|\bar{P}|)\geq
\left\lfloor\frac{3k}{2}\right\rfloor+|\bar{P}|.
$$
Hence, $(H\cup \bar{P},C)\in \mathcal{F}_{k}(|\bar{P}|+1)$ and
$(H\cup \bar{P}[S](E\cup S),K(S))\in
\mathcal{F}^{+}_{k}(|\bar{P}|+1)$. From Lemma \ref{le3}, there
exists a path $Q\subseteq E$ of order $|\bar{P}|+1$ starting from
$q$ such that $(H\cup \bar{P})[S](E\cup S)-Q$ is $k$-connected. Now
for $(H\cup \bar{P},C)\in \mathcal{F}_{k}(|\bar{P}|+1)$ by Lemma
\ref{le2}(c), then $\kappa((H\cup \bar{P})[S](E\cup S)-Q)\geq k$.
The spider
$T^{t;j+1}_{t+j+2}:=(T^{t;j}_{t+j+1}\setminus
\bar{P})\cup yq\cup Q$ rooted at $s_{0}$ has order
$|T^{t;j}_{t+j+1}|+1\leq m$ and
$G-V(T^{t;j+1}_{t+j+2})=H-Q$ is $k$-connected, a contradiction.

We can again prove it holds for $y\in \{s_{1},s_{2},\cdots,s_{t}\}$ by the above way, where $\bar{P}=P[p_{j},y)$.  Complete the proof of Lemma \ref{thm1}.

\begin{claim}\label{cl6}
There is a lager spider $T^{t+1;j}_{t+j+2}$ such that $G-T^{t+1;j}_{t+j+2}$ remains $k$-connected.
\end{claim}
\begin{proof}

Since $d_{G}(s_{0})=|G|-1$, it follows that $N(s_{0})\cap E\neq \emptyset$. Let $N(s_{0})\cap E=\{q\}$. Furthermore, from Lemma \ref{le3}, there
exists a vertex $q\subseteq E$ such that $H[S](E\cup S)-q$ is $k$-connected. Hence, $G-V(T^{t+1;j}_{t+j+2})=H-q$ is $k$-connected, a contradiction.
\end{proof}
By Claims \ref{cl5} and \ref{cl6}, we completed the proof.
\end{proof}

\section{Concluding remark}

As we all know, the higher the connectivity is, the more complex the structure of the graph is. In order to complete the proof of Mader's conjecture, we need to prove more structure theorems. According to our research experience, it needs a totally new method if one would to confirm thoroughly Mader's conjecture for the highly connected graph $G$. Naturally, maybe one could replace the minimum degree condition with some others, such as degree sum condition. In addition, one could also confirm it for some special $k$-connected graphs $G$.

Maybe one could prove it for more special trees utilizing our results. For the tree $T^{t;i,m-t-1-i}_{m}$, it means $T^{t;i,m-t-1-i}_{m}$ has $t$ legs of length one, one leg of length $i$, and one leg of length $m-t-i-1$.

\begin{problem}
Every $k+1$-connected graph $G$ with
$\delta(G)\geq\lfloor\frac{3k}{2}\rfloor+m-1$ and $\Delta(G)=|G|+1$ for positive integers
$k,m,t$, contains a spider $T^{t;i,m-t-1-i}_{m}$ such that
$$\kappa(G-V(T^{t;i,m-t-i-1}_{m}))\geq k.$$
\end{problem}

\noindent{\bf Acknowledgement} We would like to thank Dr. Shinya Fujita for the valuable suggestions and comments.
\vskip 2mm
\noindent{\bf Declaration of interests}

The authors declare that they have no known competing financial interests or personal relationships that could have appeared to influence the work reported in this paper.

\end{document}